\documentclass[12pt]{article}
\usepackage{amsfonts}

 \def\medskipamount{12pt} \def\smallskipamount{6pt}
 \def\arraystretch{1.6}

\newcounter{bitcount}
\newcommand{\bit}[1]{\addtocounter{bitcount}{1}\pagebreak[3]
\subsection{#1}\nopagebreak\setcounter{equation}{0}}

\renewcommand{\theequation}{\thesubsection .\arabic{equation}}
\renewcommand{\thesubsection}{\arabic{bitcount}}

\newcommand{\re}[1]{\mbox{\bf (\ref{#1})}}

\catcode`\@=\active
\catcode`\@=11
\def\@eqnnum{\hbox to .01pt{}\rlap{\bf \hskip -\displaywidth(\theequation)}}
\catcode`\@=12

\begin{document}


\catcode`\@=\active
\catcode`\@=11
\newcommand{\nc}{\newcommand}


\nc{\bs}[1]{ \addvspace{\medskipamount} \pagebreak[3]
\refstepcounter{equation}
\noindent {\bf (\theequation) #1.} \begin{em} \nopagebreak}

\nc{\es}{\end{em} \par \addvspace{\medskipamount} } 

\nc{\ess}{\end{em} \par}

\nc{\br}[1]{ \addvspace{\medskipamount} \pagebreak[3]
\refstepcounter{equation} 
\noindent {\bf (\theequation) #1.} \nopagebreak}

\nc{\brs}[1]{ \pagebreak[3]
\refstepcounter{equation} 
\noindent {\bf (\theequation) #1.} \nopagebreak}

\nc{\er}{\par \addvspace{\medskipamount} }


\nc{\vars}[2]
{{\mathchoice{\mb{#1}}{\mb{#1}}{\mb{#2}}{\mb{#2}}}}
\nc{\C}{\mathbb C}
\renewcommand{\H}{H}
\nc{\Q}{\mathbb H}
\nc{\R}{\mathbb R}
\nc{\Z}{\mathbb Z}
\renewcommand{\P}{\mathbb P} 

\renewcommand{\a}{{\mathfrak a}}
\renewcommand{\t}{{\mathfrak t}}


\nc{\oper}[1]{\mathop{\mathchoice{\mbox{\rm #1}}{\mbox{\rm #1}}
{\mbox{\rm \scriptsize #1}}{\mbox{\rm \tiny #1}}}\nolimits}
\nc{\Ext}{\oper{Ext}}
\nc{\Hom}{\oper{Hom}}
\nc{\Tor}{\oper{Tor}}

\nc{\Dol}{{\oper{Dol}}}
\nc{\DR}{{\oper{DR}}}
\nc{\B}{{\oper{B}}}

\nc{\h}{{\mathfrak h}}

\nc{\D}{J}

\nc{\operlim}[1]{\mathop{\mathchoice{\mbox{\rm #1}}{\mbox{\rm #1}}
{\mbox{\rm \scriptsize #1}}{\mbox{\rm \tiny #1}}}}

\nc{\lp}{\raisebox{-.1ex}{\rm\large(}}
\nc{\rp}{\raisebox{-.1ex}{\rm\large)}}

\nc{\al}{\alpha}
\nc{\be}{\beta}
\nc{\la}{\lambda}
\nc{\La}{\Lambda}
\nc{\ep}{\varepsilon}
\nc{\si}{\sigma}
\nc{\om}{\omega}
\nc{\Om}{\Omega}
\nc{\Ga}{\Gamma}
\nc{\Si}{\Sigma}


\nc{\Left}[1]{\hbox{$\left#1\vbox to
    11.5pt{}\right.\nulldelimiterspace=0pt \mathsurround=0pt$}}
\nc{\Right}[1]{\hbox{$\left.\vbox to
    11.5pt{}\right#1\nulldelimiterspace=0pt \mathsurround=0pt$}}

\nc{\updown}{\hbox{$\left\updownarrow\vbox to
    10pt{}\right.\nulldelimiterspace=0pt \mathsurround=0pt$}}


\nc{\beqas}{\begin{eqnarray*}}
\nc{\co}{{\cal O}}
\nc{\cx}{{\C^{\times}}}
\nc{\down}{\Big\downarrow}
\nc{\Down}{\left\downarrow
    \rule{0em}{8.5ex}\right.}
\nc{\downarg}[1]{{\phantom{\scriptstyle #1}\Big\downarrow
    \raisebox{.4ex}{$\scriptstyle #1$}}}
\nc{\eeqas}{\end{eqnarray*}}
\newsymbol\quadri 1003
\nc{\fp}{\mbox{     $\quadri$} \par \addvspace{\smallskipamount}}
\nc{\lrow}{\longrightarrow}
\nc{\pf}{\noindent {\em Proof}}
\nc{\sans}{\, \backslash \,}
\nc{\st}{\, | \,}

\nc{\mod}{\, / \,}

\nc{\Eorb}{E_{\oper{orb}}}

\nc{\M}{{\mathcal M}}
\nc{\N}{{\mathcal N}}

\nc{\g}{{\mathfrak g}}
\nc{\beq}{\begin{equation}}
\nc{\eeq}{\end{equation}}

\nc{\gp}[2]{{{\rm #1(}#2{\rm )}}}
\nc{\PGL}[1]{\gp{PGL}{#1}}
\nc{\GL}[1]{\gp{GL}{#1}}
\nc{\SL}[1]{\gp{SL}{#1}}
\nc{\U}[1]{\gp{U}{#1}}

\hyphenation{para-met-riz-ing sub-bundle}

\catcode`\@=12



\noindent
{\LARGE \bf Mirror symmetry, Langlands duality, }
\smallskip \\ 
{\LARGE \bf and commuting elements of Lie groups}
\medskip \\ 
{\bf Michael Thaddeus } \smallskip \\ 
Department of Mathematics, Columbia University \\
2990 Broadway, New York, N.Y. 10027
\renewcommand{\thefootnote}{}
\footnotetext{Partially supported by NSF grant DMS--9808529.}

{\small
\begin{quote}
\noindent {\em Abstract.}
By normalizing a component of the space of commuting pairs of elements
in a reductive Lie group $G$, and the corresponding space for the
Langlands dual group, we construct pairs of hyperk\"ahler orbifolds
which satisfy the conditions to be mirror partners in the sense of
Strominger-Yau-Zaslow.  The same holds true for commuting quadruples
in a compact Lie group.  The Hodge numbers of the mirror partners, or
more precisely their orbifold $E$-polynomials, are shown to agree, as
predicted by mirror symmetry.  These polynomials are explicitly
calculated when $G$ is a quotient of $\SL{n}$.
\end{quote}
}

Mirror symmetry made its first appearance in 1990 as an equivalence
between two linear sigma-models in superstring theory \cite{cls,gp}.
The targets were Calabi-Yau 3-folds, so mirror symmetry predicted that
these should come in pairs, $M$ and $\hat M$, satisfying $h^{p,q}(M) =
H^{p,3-q}(\hat M)$.

Although many examples were known, the physics did not immediately
provide any general construction of a mathematical nature for the
mirror.  Since then, however, two mathematical constructions have
emerged: that of Batyrev \cite{b} and Batyrev-Borisov \cite{bb},
generalizing the original idea of Greene-Plesser \cite{gp}, and that of
Strominger-Yau-Zaslow \cite{syz} with which this paper is concerned. 

Of the two, Batyrev's construction has the advantage of being precise,
and more amenable to explicit calculations.  One can prove, for example,
that the Hodge numbers of the Batyrev mirror satisfy the desired
relationship.  On the other hand, it is deeply rooted in toric
geometry.  This has led skeptics to suggest that mirror symmetry is
an intrinsically toric phenomenon, despite work \cite{bckv,einar}
extending Batyrev's point of view some ways beyond the toric setting.

The construction proposed by Strominger-Yau-Zaslow in 1996 has quite a
different flavor.  It is directly inspired by a physical duality,
the so-called $T$-duality between sigma-models whose targets are dual
tori.  Remarkably, although it is supposed to transform one
projective variety into another, the construction is not algebraic, or
even K\"ahler, in nature.  Rather, it is symplectic: one must
find a foliation of $M$ by special Lagrangian tori, and replace
each torus with its dual.  

This bold idea has already led to some interesting work on the existence
of families of Lagrangian tori in Calabi-Yau 3-folds \cite{g1,g2,r},
which is essentially a problem in symplectic topology. But it is not
yet sufficiently advanced that the mirror can be constructed in any
precise sense, nor any of its invariants computed beyond the Euler
characteristic.  It is not even known how to construct families of
tori which are special Lagrangian (as opposed to just Lagrangian).
And the further questions of what complex structure to place on the
dual family, and how to deal with singular fibers, remain mysterious.

This paper will study a case where these formidable difficulties can
be completely circumvented.  The Strominger-Yau-Zaslow construction
works like a charm, and the equality of the Hodge numbers can be
completely verified.  Furthermore, there is a suggestive relationship
with the Langlands duality between reductive Lie groups.  On the other
hand, this case is in some ways quite distant from the original case
of Calabi-Yau 3-folds.  The spaces in question are singular and often
non-compact.  Moreover, their complex dimension is even and typically
rather large.

To find the special Lagrangian tori, we exploit the existence of a
{\em hyperk\"ahler} structure.  This involves three complex structures
on $M$ which satisfy the commutation relations of the imaginary
quaternions.  It is easy to show that submanifolds which are
holomorphic Lagrangian with respect to one complex structure are
special Lagrangian with respect to another.  Hence we can study
special Lagrangian fibrations without leaving the realm of algebraic
geometry.  If $M$ is compact and hyperk\"ahler, its Hodge numbers
satisfy $h^{p,q}(M) = h^{n-p,q}(M)$, where $n$ is the complex
dimension of $M$.  So mirror symmetry leads us to expect that
$h^{p,q}(M) = h^{p,q}(\hat M)$.  We will see that this is indeed the
case.  It even remains true for the non-compact hyperk\"ahler examples
we shall encounter.

The metric inducing this hyperk\"ahler structure is not only
hyperk\"ahler, it is flat!  So the Riemannian geometry of the
situation is not very interesting.  What makes things non-trivial is
the presence of orbifold singularities.  In a few of our cases, these
can be resolved crepantly \cite{bdl}, leading to Beauville's examples
\cite{beau} of compact hyperk\"ahler manifolds, which do carry
non-flat metrics.  But in most cases, there is no crepant resolution.
What we actually evaluate, therefore, are not actual Hodge numbers,
but rather orbifold Hodge numbers in the sense of Vafa \cite{v},
Zaslow \cite{zas}, and Batyrev-Dais \cite{bd}.

It is worth emphasizing that, although many aspects of mirror symmetry
remain highly speculative, the results on orbifold Hodge numbers
in this paper are mathematically rigorous.  Their motivation and
interpretation are, of course, more open-ended.  To clarify this
distinction, the paper has been divided into two parts.  The first,
comprising sections 1 to 3, consists of background and motivation.  The
second, comprising sections 4 to 7, contains precise mathematical
statements and proofs.

Section 1 gives the necessary background on hyperk\"ahler manifolds
and special Lagrangian tori.  Sections 2 and 3 introduce the two main
classes of examples we shall study.  The first is a component of the
space of conjugacy classes of commuting pairs of elements in a complex
reductive group $G$.  It can be identified, thanks to the work of
Hitchin, with the space of Higgs bundles on an elliptic curve.  The
second, likewise, is a component of the space of conjugacy classes of
commuting {\em quadruples\/} of elements of a {\em compact\/} Lie
group $K$.  In each case, the mirror is the same kind of space, only
with the group replaced by its Langlands dual.

The next three sections consist of a rigorous formulation and proof of
the equality of Hodge numbers predicted by mirror symmetry for these
spaces.  One first has to normalize the spaces, which is carried out
in section 4.  Then, since they are not compact or smooth, one needs
to define the Hodge numbers judiciously: the suitable definitions, of
the so-called orbifold $E$-polynomials, are recalled in section 5.
Finally, the main theorem, showing that the orbifold $E$-polynomials of
the putative mirror partners agree, is proved in section 6.

Section 7 consists of a more or less explicit evaluation of the
orbifold $E$-polynomials for a group of the form $\SL{n}/\Z_m$.  Section
8 contains some concluding remarks and suggestions for further
research. \smallskip

\noindent {\em Notation.} A few conventions should be
mentioned.  First, all varieties are over the complex numbers.
Second, $\otimes$ means $\otimes_\Z$ unless otherwise
stated.  Finally, the expression ``Killing form'' is used loosely to
mean any nondegenerate symmetric bilinear form on a complex reductive
Lie algebra $\g$ which is compatible with the splitting $\g =
{\mathfrak z}(\g) \oplus [\g,\g]$ and restricts to the Killing form, in
the usual sense, on the second factor. \smallskip

\noindent {\em Acknowledgements.} This work was inspired by several
conversations with Jim Bryan.  I wish to thank him for his generosity
with his ideas, and for his continuing interest and
advice.  I am also grateful to Lev Borisov, Ron Donagi, Robert
Friedman, Dennis Gaitsgory, Paul Gunnells, John Morgan, Carlos
Simpson, and the referee for helpful advice.

\bit{Special Lagrangian tori on hyperk\"ahler manifolds}

Let $M$ be a Calabi-Yau manifold of real dimension $2n$, with K\"ahler
form $\om$ and holomorphic $n$-form $\Om = \Om_1 + i \Om_2$.  A real
submanifold $L$ of $M$ is said to be {\em Lagrangian} if $\dim L = n$
and $\om|_L = 0$, and {\em special Lagrangian} if $\Om_2|_L = 0$ as
well.

To describe the mirror $\hat M$, Strominger-Yau-Zaslow ask for a real
$n$-manifold $N$ and a map $M \to N$ whose fiber over a general point
$x \in N$ is a special Lagrangian $n$-torus $L_x$.  The mirror should
be the total space of the family over $N$ whose fiber over $x$ is the
dual torus $\hat L_x = \Hom (\pi_1(L_x), \U{1})$.  At least when the
map $M \to N$ has a section, the mirror of $\hat M$ can be identified
with the original $M$ by double duality.  It is not clear, though, how
to dualize the singular fibers of the map, or how to define a complex
structure on $\hat M$.

Now suppose that $M$ is a {\em hyperk\"ahler} manifold.  This means
that it has a metric which is K\"ahler with respect to three complex
structures $J_1,J_2,J_3:TM \to TM$ satisfying the commutation
relations of the imaginary quaternions.  Let $\om_1, \om_2, \om_3$ be
the three K\"ahler forms.  Then it is not hard to show that $\om_2 + i
\om_3$ defines a complex symplectic form on $M$, which is holomorphic
when $M$ is equipped with the complex structure $J_1$.  In particular,
$M$ must have real dimension $4m$ for some $m$, and $\Om = (\om_2 +
i \om_3)^m$ is a nowhere vanishing section of $K_M$.  Therefore $M$
may be regarded as a (possibly non-compact) Calabi-Yau manifold.
Also, since $M$ is complex symplectic, $TM \cong T^*M$
holomorphically, so $\Om^p M \cong \Om^{2m-p} M$, and hence in the
compact case $h^{p,q}(M) = h^{2m-p,q}(M)$.  

We will occasionally refer to the identity $M \to M$, regarded as a
non-holomorphic map between different complex structures, as {\em
hyperk\"ahler rotation}.

Let $L$ be a real submanifold of $M$ which is complex Lagrangian with
respect to $\om_2 + i \om_3$.  Then $L$ is actually a complex
submanifold with respect to $J_1$ \cite{hit99}.  Furthermore, since
$\om_2|_L = \om_3|_L =0$, $L$ is real Lagrangian with respect to
$\om_3$, and $(\om_1+i\om_2)^m$ restricts to $L$ as $\om_1^m$, so that
its imaginary part vanishes.  Hence {\em $L$ is special Lagrangian
  with respect to $J_3$}.

This fact will be of key importance for us.  It means that, if we seek
special Lagrangian torus fibrations on a hyperk\"ahler manifold, we
need look no further than {\em holomorphic} Lagrangian torus
fibrations in a different complex structure.  In particular, the hard
analysis usually appearing in the search for special Lagrangian
submanifolds can be completely avoided: we can work in a purely
algebraic setting.

Even better, there is an obvious source of families of holomorphic
Lagrangian torus fibrations on hyperk\"ahler manifolds.  Liouville's
theorem from Hamiltonian mechanics assures us that, if we have $m$
Poisson-commuting holomorphic functions on $M$ whose derivatives are
generically independent, then the map $M \to \C^m$ they define will
have exactly the desired property, provided at least that it is
proper.  If so, we have what is called an algebraically completely
integrable Hamiltonian system.

There is a celebrated system of this kind, the so-called {\em Hitchin
system}. It lives on the moduli space of Higgs bundles on an
algebraic curve.  In this paper we shall assume that the curve is
elliptic.  However, curves of genus $>1$ will also furnish many
interesting examples of mirror pairs.  This is the subject of a
forthcoming paper of T.~Hausel and the author \cite{ht}. 

The spaces we shall consider here will actually be orbifolds, not
manifolds.  Indeed, their singularities will be exactly what makes the
computation of the Hodge numbers non-trivial.  However, they will be
global quotients of hyperk\"ahler manifolds --- indeed, of extremely
simple ones.

The torus fibrations coming from algebraically integrable systems are
such convenient sources of special Lagrangian tori that it is a little
mysterious why they have never been studied before in connection with
Strominger-Yau-Zaslow.  One reason might be the work of Verbitsky
\cite{v} on mirror symmetry and hyperk\"ahler manifolds.  Verbitsky
showed that a compact hyperk\"ahler manifold is mirror to itself in
quite a strong sense: not only are the Hodge numbers self-mirror,
which is trivial, but the Yukawa coupling also corresponds to the
quantum product (which equals the ordinary cup product in the
hyperk\"ahler case), and even the local variations of these structures
correspond.  This seems to make the search for a mirror into a
triviality.  However, our examples show, without contradicting
Verbitsky's results, that the Strominger-Yau-Zaslow mirror of a
compact hyperk\"ahler orbifold, such as the fiber of the sum map
$\oper{Sym}^m A \to A$ for an abelian surface $A$, can be a different
space: see Proposition \re{exam} below.

\bit{Higgs bundles on an elliptic curve}

Here, then, is the situation we want to consider.  Let $G$ be a
connected complex reductive algebraic group, let $\H \subset G$ be a
maximal torus, let $W = N(\H)/\H$ be the Weyl group, and let $\La
\subset \h$ be the coweight lattice.  Also let $\hat \H = \Hom(\H,
\cx)$ and $\hat \La = \Hom(\La, \Z)$.  There is in some sense
canonically associated to $G$ a reductive group $\hat G$, the {\em
Langlands dual}, with maximal torus $\hat \H$ and coweight lattice
$\hat \La$.

Let $C$ be an elliptic curve.  We will examine the moduli spaces
studied by Simpson: $M_\Dol(G)$, the moduli space of topologically
trivial semistable Higgs $G$-bundles on $C$; $M_\DR(G)$, the moduli
space of topologically trivial local systems on $C$ with structure
group $G$; and $M_\B(G)$, the identity component of the moduli space
of homomorphisms $\pi_1(C) \to G$ modulo conjugacy.  The subscripts
stand for Dolbeault, de Rham, and Betti respectively, and are
Simpson's notation.

We refer to Simpson \cite{s2,s} for the precise definitions of these
spaces, including the correct notions of stability and equivalence.
Suffice it to say that a {\em Higgs $G$-bundle} is a pair $(E,\phi)$,
where $E$ is a holomorphic principal $G$-bundle, and $\phi \in
H^0(C,\oper{ad} E \otimes K_C)$.  This is particularly simple in the
present case, since $K_C$ is trivial.  A {\em local system} is a
holomorphic principal $G$-bundle with an integrable holomorphic
connection.  A local system is determined up to isomorphism by its
holonomy; this induces an isomorphism $M_\DR(G) \cong M_\B(G)$, but
only analytically, since the universal cover of $C$ must be used to
construct the isomorphism.

Note also that since $\pi_1(C) = \Z \oplus \Z$, $M_\B(G)$ is simply
the identity component of the space of commuting pairs of elements of
$G$, modulo conjugacy.

Simpson \cite{s2} showed that there is a natural homeomorphism between
$M_\Dol(G)$ and $M_\DR(G)$.  It is not a biholomorphism on the smooth
locus; rather, Hitchin \cite{hit87} showed that the two complex
structures $J_1, J_2$ coming from $M_\Dol(G)$ and $M_\DR(G)$ form part
of a hyperk\"ahler structure.

There is a dominant proper morphism, the {\em Hitchin map}, which
takes $M_\DR(G)$ to a vector space, and is defined by evaluating a set
of generators for the invariant polynomials of $\g$ on $\phi$.  The
map in question therefore takes values in $\g/G$, which according
to a classic theorem of Chevalley \cite[4.9]{vara} is a vector space
and is isomorphic to $\h/W$.

The fibers of the Hitchin map have been much studied and are known to
be generically abelian varieties \cite{duke}.  In our case, the fiber
over a regular element of $\h/W$ can be identified as $C \otimes \La$.
This can, for example, be extracted from the work of Donagi-Gaitsgory
\cite{dg} or Faltings \cite{falt}: according to their recipe one gets
$H^1(C,\co(\H))$, where $\co(\H)$ denotes the sheaf of germs of
holomorphic maps to $\H$.  Since $\H = \cx \otimes \La$, this is
naturally isomorphic to $H^1(C,\co^\times) \otimes \La$; but since we
are dealing only with topologically trivial bundles, we must take only
the identity component, which is $C \otimes \La$.

How are the moduli spaces for $G$ and $\hat G$ related?  Exactly as
one would hope for the Strominger-Yau-Zaslow construction to apply!
Since $W$ acts orthogonally with respect to the Killing form, $\h
\cong \hat \h$ as representations of $W$.  So $\h/W \cong \hat \h/W$:
the Hitchin maps have isomorphic targets.  On the other hand, $C
\otimes \La = \oper{Pic}^0(C \otimes \hat \La)$: the generic fibers of
the two Hitchin maps are dual.  Furthermore, the Higgs bundles whose
underlying bundles are trivial furnish a section of each Hitchin map.

Hence $M_\DR(G)$ and $M_\DR(\hat G)$ satisfy the requirements to be
Strominger-Yau-Zaslow mirror partners.  We might conjecture,
therefore, that their Hodge numbers are equal.  But what is meant by
the Hodge numbers of varieties which are not compact or smooth?  The
right answer seems to be the {\em stringy $E$-polynomial} of
Batyrev-Dais \cite{bd}.  Our chief goal, attained in \S6, will be
to show that our moduli spaces have the same stringy $E$-polynomials
as their putative mirror partners.

One technicality stands in the way, however.  The minimum requirements
for varieties to have well-defined stringy $E$-polynomials
are not altogether clear, but they certainly must be {\em normal}.  It
is by no means obvious that the moduli spaces discussed above meet
this minimum requirement.

The solution is to replace every space in the story by its
normalization.  This motivates \S4, in which it is shown that the
normalizations of $M_\Dol(G)$, $M_\DR(G)$, and $M_\B(G)$ are
$(T^*C~\otimes~\La)/W$, $(H\times H)/W$, and $(H\times H)/W$
respectively.  The cotangent bundle $T^*C = C \times \C$ has an
obvious hyperk\"ahler structure, which induces one on $T^*C \otimes
\La$, so $(T^*C~\otimes~\La)/W = ((C \otimes \La) \times \h)/W$ is a
hyperk\"ahler orbifold.  The second complex structure on $T^*C$ is
$\cx \times \cx$, so the second complex structure on $(T^*C \otimes
\La)/W$ is $((\cx \times \cx) \otimes \La)/W = (H \times H)/W$.  The
Hitchin map lifts to the projection $((C \otimes \La) \times \h)/W \to
\h/W$, whose generic fiber is clearly $C \otimes \La$.  The mirror
transformation can therefore be summarized as follows.  The vertical
arrows denote hyperk\"ahler rotation, and the diagonal arrows are the
normalized Hitchin maps.
$$\def\arraystretch{2.2}
\begin{array}{ccc}
\displaystyle \frac{\H \times \H}{W} 
&& \displaystyle \frac{\hat \H \times \hat \H}{W} \\
\updown && \updown \\
\displaystyle \frac{T^*C \otimes \La}{W} 
&& \displaystyle \frac{T^*C \otimes \hat \La}{W} \\
\phantom{XXX}\searrow && \swarrow \phantom{XXX} \\
&\displaystyle \frac{\h}{W} &
\end{array}$$ 
We hope, and indeed will show in \S6, that the stringy
$E$-polynomials of the spaces in the top row agree.

\bit{Principal bundles on an abelian surface}

Discussions with Jim Bryan in October 1998 revealed that there is an
analogue of the above situation where all the spaces involved are
compact.  They are, in fact, precisely the spaces studied in the
recent work of Bryan-Donagi-Leung \cite{bdl}.

Let $A = \U{1}^4$ be a 4-torus with a hyperk\"ahler structure
compatible with the flat metric.  This is easily found: just choose
any linear isomorphism from ${\mathfrak u}(1)^4$ to the quaternions.
For convenience, we will suppose that in one complex structure, $A$
splits as a product $C \times D$ of two elliptic curves.  However,
many of our considerations would extend to the more general case of an
abelian surface which is an extension of one elliptic curve by
another.

Consider the moduli space of topologically trivial semistable
principal $G$-bundles on $A$.  Let $M(G)$ be the connected component
containing the trivial bundle.  The correspondence theorem relating
holomorphic bundles to Hermitian-Einstein connections \cite{don,uy}
has been generalized to arbitrary compact structure groups by
Ramanathan-Subramanian \cite{rs}.  So if $K$ is a
compact Lie group whose complexification is $G$, then every
equivalence class of topologically trivial semistable $G$-bundles has
a representative carrying a flat $K$-connection, unique up to
conjugacy.

Hence $M(G)$ is homeomorphic to a connected component of the space of
homomorphisms $\pi_1(A) \to K$ modulo conjugacy, that is, the space of
commuting {\em quadruples} of elements of $K$, modulo conjugacy.  As
shown by Borel-Friedman-Morgan \cite[2.3.2]{bfm}, the connected
component of the identity consists precisely of those quadruples which
all belong to a common maximal torus.

Since $A$ has several complex structures, this endows the moduli space
$M(G)$ with several complex structures.  One would expect that they
give rise to a hyperk\"ahler structure on the smooth locus, but this
result does not seem to appear in the literature except for $G =
\GL{n,\C}$ \cite{m2}.  However, this lacuna will not be serious as
the normalization of $M(G)$ will manifestly be a hyperk\"ahler
orbifold.  In fact, it is shown in \S4 that it is exactly the quotient
$(\hat A \otimes \La)/W$.  In the present case this can be identified
with  $(A \otimes \La)/W$, since $C \times D$ is self-dual and
hence so is $A$ in the other K\"ahler structures.

When $A$ has the complex structure of $C \times D$, then restriction
to $C$ times a point defines a morphism from $M(G)$ to the moduli space
of semistable $G$-bundles on $C$.  (That semistability is preserved
follows from the correspondence theorem.)  The fibers of this morphism
are not as well-studied as those of the Hitchin map.  But the
situation is clarified by passing to the normalization.  The space of
topologically trivial semistable $G$-bundles on an elliptic curve $C$
is shown by Friedman-Morgan \cite{fm} and Laszlo \cite{laszlo} to be
$(C \otimes \La)/W$.  (When $G$ is simply connected, this is the space
shown by Looijenga \cite{loo} and Bernshtein-Shvartsman \cite{bs} to
be a weighted projective space.)  The normalization of the restriction
map is easily seen to be just the map $(A \otimes \La)/W \to (C
\otimes \La)/W$ induced by the projection $A \to C$.

The generic fiber is clearly the abelian variety $D \otimes \La$.
Moreover, there is an obvious candidate for the dual fibration, namely
the quotient of $(D \otimes \hat \La) \times (C \otimes \La)$ by
$W$.  Since this torus has tangent space naturally isomorphic to that
of $A \otimes \La$, it too carries a natural hyperk\"ahler structure.
Denote ${\mathcal B}$ the space obtained by hyperk\"ahler rotation.
It does not split naturally, but it nevertheless carries an action of
$W$, and ${\mathcal B}/W$ is the mirror of $(A \otimes \La)/W$.

The story so far is summarized by the left-hand side of the diagram
below, with the vertical arrows denoting hyperk\"ahler rotation as
before.  But the roles of $\La$ and $\hat \La$ are completely
interchangeable, so the diagram has a right-hand side as well.  For
that matter, the diagram could continue with $\hat {\mathcal B}/W$,
to make a closed chain with four links. 
$$\def\arraystretch{2.2}
\begin{array}{ccccc}
\displaystyle \frac{A \otimes \La}{W} 
&& \displaystyle \frac{\mathcal B}{W} 
&& \displaystyle \frac{A \otimes \hat \La}{W} \\
\updown && \updown && \updown\\
\displaystyle \frac{(C \times D) \otimes \La}{W} 
&& \displaystyle \frac{(D \otimes \hat \La) \times (C \otimes \La)}{W} 
&& \displaystyle \frac{(C \times D) \otimes \hat \La}{W} \\
\phantom{XXXX}\searrow 
&& \swarrow \phantom{XXXXXX} \searrow 
&& \swarrow \phantom{XXXX} \\
&\displaystyle \frac{C \otimes \La}{W} 
&& \displaystyle \frac{D \otimes \hat \La}{W} &
\end{array}$$ 
Thus each space in the top row has {\em two} special Lagrangian torus
fibrations, and hence two mirrors! Under the circumstances, it is
tempting to guess that perhaps these two mirrors are isomorphic.
However, this is generally not the case.  Indeed, we shall see in
Proposition \re{exam} that for $G = \SL{3}$, $(A \otimes \La)/W$ has a
crepant resolution while $(A \otimes \hat \La)/W$ does not.

Once again, we will be able to verify the prediction of mirror
symmetry that the stringy $E$-polynomials of these spaces
will be the same.

\bit{The normalizations of the moduli spaces}

As in the previous two sections, let $G$ be a complex connected
reductive group, $\La$ its coweight lattice, $W$ the Weyl group, and
let $C$ be an elliptic curve.  We will describe the normalizations of
Simpson's moduli spaces over $C$ as follows.

\bs{Theorem}
\label{norm}
The normalizations of $M_\Dol(G)$, $M_\DR(G)$, and $M_\B(G)$ are
naturally isomorphic to $(T^*C \otimes \La)/W$, $(\D \otimes \La)/W$, and
$(H \times H)/W$ respectively.  Here $\D$ denotes the unique algebraic
group which is a non-trivial extension of $C$ by the affine line.
\es

To construct the morphisms that will turn out to be the
normalizations, consider first the case $G = \cx$.  Clearly
$M_\Dol(\cx) = T^*C = C \times \C$ and $M_\B(\cx) = \cx \times \cx$.
Moreover, there is an algebraic group epimorphism from $M_\DR(\cx)$ to
the moduli space of topologically trivial holomorphic $\cx$-bundles on
$C$, which of course is isomorphic to $C$.  The kernel consists of the
holomorphic connections on the trivial bundle, which can be identified
with $H^1(C,\co) \cong \C$.  The extension is not split, since
$M_\DR(\cx)$ is analytically isomorphic to $M_\B(\cx) = \cx \times
\cx$ which cannot contain a complete curve.  But there is well known
to be only a single non-trivial extension $\D$ of $C$ by $\C$: see for
example Serre \cite[VII \S3]{serre}.

Now for general connected reductive $G$, let $H= \cx \otimes \La$ be a
Cartan subgroup.  Then $M_\Dol(H) = T^*C \otimes \La$, $M_\DR(H) = \D
\otimes \La$, and $M_\B(H) = (\cx \times \cx) \otimes \La = H \times
H$.  Extension of the structure group induces a morphism $M_\Dol(H)
\to M_\Dol(G)$, and likewise for the de Rham and Betti spaces.  The
Weyl group acts by outer automorphisms on $H$ that extend to inner
automorphisms on $G$, so it acts on $M_\Dol(H)$ via its action on
$\La$, but trivially on $M_\Dol(G)$.  The aforementioned morphism
therefore descends to a morphism $\rho_\Dol: (T^*C \otimes \La)/W \to
M_\Dol(G)$.  Similarly, there are morphisms $\rho_\DR$ and $\rho_B$.

Note that there is a commutative diagram of continuous maps
$$
\def\arraystretch{2.2}
\begin{array}{ccccc}
\displaystyle \frac{T^*C \otimes \La}{W} & \longleftrightarrow &
\displaystyle \frac{\D \otimes \La}{W} & \longleftrightarrow &
\displaystyle \frac{H \times H}{W} \\
\downarg{\rho_\Dol} && \downarg{\rho_\DR} && \downarg{\rho_\B} \\
M_\Dol(G) & \longleftrightarrow &
M_\DR(G) & \longleftrightarrow &
M_\B(G).
\end{array}$$

\bs{Proposition}
\label{finite}
The morphisms $\rho_\Dol$, $\rho_\DR$ and $\rho_\B$ are finite
bijections. 
\es

\pf.  Being finite is equivalent to being proper with finite fibers.
Everything is therefore a topological statement.  Since the
rows in the diagram above consist of homeomorphisms, it suffices to
prove the proposition for one of the three moduli spaces.  We select
the Betti space.

Both $(H \times H)/W$ and $M_\B(G)$ have natural morphisms to $\H/W
\times \H/W$.  The projection $f: (\H \times \H)/W \to \H/W \times \H/W$ is
the obvious one, and the morphism $g:\mu^{-1}(I)/G \to \H/W \times \H/W$
is equally obvious once $\H/W$ is identified \cite[6.4]{stein} with the
quotient $G/G$, where $G$ acts on itself by conjugation.  These
morphisms satisfy $f = g \circ \rho_\B$.  And $f$ is finite: the
coordinate ring of $\H \times \H$ is finitely generated over that of
$\H/W \times \H/W$, so the same is true of its submodule, the coordinate
ring of $(\H \times \H)/W$.  Hence $\rho$ is a proper morphism \cite[II
4.8(e)]{h} of affine varieties, and hence finite \cite[II Ex.\ 
4.6]{h}.

The work of Richardson \cite{rich} implies that if $G$ is reductive,
connected and simply connected, then the set of commuting pairs in $G
\times G$ is irreducible, and hence $M_\B(G)$ is irreducible.  It
follows that $M_\B(G)$ is still irreducible even if $G$ is not simply
connected.  For it was defined to be the connected component of the
space of commuting pairs containing $(I,I)$.  This is surjected on by
$Z_0(G)^2 \times M_\B(\tilde G)$, where $Z_0(G)$ is the identity
component of the center, and $\tilde G$ is the universal cover of the
commutator subgroup.

Hence to show that the proper map $\rho_\B$ is surjective, it suffices
to show that it is dominant.  Let $h_1,h_2$ be any commuting regular
semisimple elements of a Cartan subgroup $\H \subset G$.  The kernels
of $\oper{ad}_{h_1}$ and $\oper{ad}_{h_2}$ on $\mathfrak g$ are both
precisely $\mathfrak h$.  Hence there exists a neighborhood of
$(h_1,h_2) \in M_\B(G)$, in the complex topology, such that for
$(g_1,g_2)$ in this neighborhood, $\ker \oper{ad}_{g_1}$ and $\ker
\oper{ad}_{g_2}$ coincide and have dimension equal to the rank of $G$.
We may also assume that $g_1$ and $g_2$ remain regular semisimple
elements in this neighborhood, so their centralizers remain tori.  The
Lie algebras of these centralizers are both $\ker \oper{ad}_{g_1}$, so
they are both the same maximal torus and hence $g_1$ and $g_2$ belong
to the same Cartan subgroup.  Hence $(g_1,g_2)$ is in the image of
$\rho_\B$.  Therefore $\rho_\B$ surjects onto a neighborhood in the
complex topology, and hence is dominant.

To show $\rho_\B$ is injective, we follow an argument of Borel
\cite{borel}.  Suppose that $\rho_\B(h_1,h_2) = \rho_\B(h'_1,h'_2)$.
This means that $(h_1,h_2) \in \H \times \H$ is conjugate to
$(h'_1,h'_2) \in H \times H$ by some $g \in G$.  The centralizer of
$h_1$ and $h_2$, say $Z(h_1,h_2) \subset G$, is a reductive subgroup.
This follows, for example, from 26.2A of Humphreys \cite{humph}, since
the proof there is valid not only for a subtorus, but for any subset.
Then $H$ and $gHg^{-1}$ are maximal tori in $Z(h_1,h_2)$, so
$gHg^{-1}$ is conjugate to $H$ by some $g' \in Z(h_1,h_2)$.  Then $g'g
\in N(H)$ conjugates $(h_1,h_2)$ to $(h'_1,h'_2)$, so they represent
the same point in $(H \times H)/W$. \fp

\noindent {\em Proof of Theorem \re{norm}}.
First, note that, as varieties surjected on by irreducible varieties, all
three moduli spaces are irreducible.

For any $x$ in the dense set of the proposition, the fiber
$\rho_\B^{-1}(x)$ has a single closed point.  Since $\rho_\B$ is
finite, the locus where its fibers are reduced is open, and, of
course, it contains the generic point.  Hence over a nonempty open
set, the fiber of $\rho$ is a single reduced point.  Therefore
$\rho_\B$ is birational.  As the quotient of a smooth variety by a
finite group, $(H \times H)/W$ is certainly normal, so $\rho_\B$ lifts
to a birational finite morphism $\tilde \rho_\B: (H \times H)/W \to
\tilde M_\B(G)$ on the normalizations.  This is an isomorphism by
Zariski's Main Theorem.

The same argument applies to the Dolbeault and de Rham spaces.  \fp

Exactly the same methods can be used to prove a similar theorem on the
space $M(G)$ of $G$-bundles on an abelian surface $A$ which was
discussed in \S3.

\bs{Theorem}
\label{q}
The normalization of $M(G)$ is naturally isomorphic to $(\hat A \otimes
\La)/W$.  
\es

The proof is parallel to that of Theorem \re{norm}.  One defines a
morphism $\rho: (\hat A \otimes \La)/W \to M(G)$, and it then suffices
to prove the following analogue of \re{finite}.  Finiteness follows
automatically since the domain is compact.

\bs{Proposition}
The morphism $\rho$ is a bijection.
\es

\pf.  This can be deduced conveniently using the
correspondence with flat connections.  So let $K$ be a compact Lie
group whose complexification is $G$, $T \subset K$ the maximal torus
whose complexification is the Cartan subgroup $H$.
Then $M(G)$ gets identified with the identity component of the space
of commuting quadruples in $K$ modulo conjugacy, and $M(H) = \hat A
\otimes \La$ gets identified with $T^4$.  

Surjectivity is a consequence of the aforementioned result of
Borel-Friedman-Morgan \cite[2.3.2]{bfm}.  For injectivity, we again
follow the argument of Borel.  If two quadruples of elements of $T$
have the same image in $M(G)$, then there exists $g \in K$ conjugating
one quadruple to another.  Let $Z \subset K$ be the centralizer of the
first quadruple.  This is compact and contains $T$ and $gTg^{-1}$ as
maximal tori, so there exists $g' \in T$ conjugating $gTg^{-1}$ to
$T$.  Then $g'g \in N(T)$ and conjugates one quadruple to another, so
the two quadruples represent the same point in $T^4/W$. \fp

To show that our mirror transformation is non-trivial, here is an
example of a space that certainly differs from its mirror.

\bs{Proposition} 
\label{exam}
The normalized moduli spaces of Theorems \re{norm} and \re{q} have
crepant resolutions for $G = \SL{3}$, but not for $G=\PGL{3}$.  \es

\pf. Let $A$ be any 2-dimensional connected abelian algebraic group,
and let $G = \SL{3}$.  The Weyl group is the symmetric group
$S_3$, and the torus $A \otimes \La$ can be regarded as the kernel of
the sum map $A^3 \to A$.  The quotient $(A \otimes \La)/S_3$ is
therefore the fiber of the sum map $\oper{Sym}^3 A \to A$, where
$\oper{Sym}^3 A$ is the symmetric product $A^3/S_3$.  As explained by
Beauville \cite{beau}, this has a crepant resolution, namely the fiber
of the corresponding map $\oper{Hilb}^3 A \to A$, where $\oper{Hilb}^3
A$ is the Hilbert scheme parametrizing subschemes of dimension 0 and
length 3.

The torus $A \otimes \hat \La$ is, of course, isomorphic to $A \otimes
\La$.  However, the action of $S_3$ on $A \otimes \hat \La$ is
different.  In terms of a basis for $\hat \La$, which splits $A
\otimes \hat \La$ as $A \times A$, the $S_3$-action is generated by
the elements $(x,y) \mapsto (-x-y,x)$ of order 3 and $(x,y) \mapsto
(-y,-x)$ of order 2.  So for any nonzero $x \in A$ with $3x=0$, $(x,x)
\in A^2$ has stabilizer $\Z_3$.  Therefore the singularity of $(A
\otimes \hat \La)/S_3$ at this point is analytically isomorphic to the
quotient $\C^4 / \Z_3$, where $\xi = e^{2 \pi i/3}$ acts by $(w,x,y,z)
\mapsto (\xi w, \xi^2 x, \xi y, \xi^2 z)$.  This has no crepant
resolution \cite[5.4]{reid}.  Hence neither does $(A \otimes \hat
\La)/S_3$: the smoothness and discrepancies of resolutions are
analytical invariants, since they can be computed on the completed
local rings.  \fp

\bit{Review of orbifold \boldmath $E$-polynomials}

Let $W$ be a finite group acting on a smooth projective variety $X$.
Then for each $p,q \geq 0$, $W$ acts on $H^{p,q}(X)$.  Denote the
character of this representation by
$h^{p,q}_W(X):W \to \C$.  One should think of the
character as a sort of equivariant Betti number. Its average value 
$$\bar h^{p,q}_W(X) = \frac{1}{|W|} \sum_{w \in W}h^{p,q}_W(X)(w)$$
is the dimension of the $W$-invariant part of $H^{p,q}(X)$.

For any $w \in W$, denote $X^w$ the
fixed-point set of $w$; this is a smooth subvariety.  Suppose 
that the action of each $w \in W$ on $K_X|_{X^w}$ is trivial.  Then Vafa
\cite{vafa} and Zaslow \cite{zas} define {\em orbifold Hodge
numbers\/} associated to the quotient orbifold $X/W$ as follows.
For any $x \in
X^w$, $w$ acts on $T_x X$ with weights
$$e^{2 \pi i \al_1}, \dots, e^{2 \pi i \al_n},$$
where $\al_1, \dots,
\al_n \in [0,1)$.  Define the {\em fermionic shift} $F(w)$ to be
$\al_1 + \cdots + \al_n$.  This is a locally constant function on
$X^w$, and it is integer-valued due to the assumption on $K_X|_{X^w}$.
To simplify the notation, assume that it is constant, as will be true
in the examples we study.  Then define the {\em orbifold Hodge
numbers} as
$$h^{p,q}_{\oper{orb}}(X/W) 
= \sum_{ \{ w \} }\bar h^{p-F(w),q-F(w)}_{C(w)}(X^w),$$
where $C(w)$ denotes the centralizer of $w$, and the sum runs over the
conjugacy classes of $W$.
 
Since some of our moduli spaces are non-compact, we wish to generalize
these notions to the case where $X$ is a smooth {\em quasi-projective}
variety.  The cohomology of $X$ then no longer carries a pure Hodge
structure, but Deligne \cite{del2,del3} constructs a canonical {\em
mixed Hodge structure\/} on the compactly supported cohomology,
$H^k_{\oper{cpt}}(X)$.  That is, there 
are two canonical filtrations on $H^k_{\oper{cpt}}(X)$, so that we may
write the Betti number as a sum 
$h^k_{\oper{cpt}} = \sum_{p,q} h^{p,q;k}$, where
$h^{p,q;k}$ are the dimensions of the quotients $H^{p,q;k}$ associated to these
filtrations \cite[2.3.7]{del2}.  In the smooth projective case,
$h^{p,q;k} = 0$ unless $p+q=k$.  If we define
$$e^{p,q} = \sum_{k \geq 0} (-1)^k h^{p,q;k},$$
then the so-called {\em $E$-polynomial}
$$E(X) = \sum_{p,q}e^{p,q}(X)u^p v^q$$
enjoys some remarkable properties.  Specifically, it is additive for
disjoint unions of locally closed subvarieties, and multiplicative for
Zariski locally trivial fibrations.

The action of a finite group $W$ on $X$ preserves the mixed Hodge
structure, since it is canonical.  Hence the spaces $H^{p,q;k}(X)$ are
representations of $W$.  As before, let $h^{p,q;k}_W(X)$ be their
characters, and let $\bar h^{p,q;k}_W(X)$ be the average values of these
characters.  Let $e^{p,q}_W$ and $E_W(X)$ be the $W$-equivariant
versions of the expressions above, and let $\bar E_W(X)$ be the same
as $E_W(X)$, but with $h$ replaced by $\bar h$. 
Note that $E_W(X \times Y) = E_W(X) E_W(Y)$, but that the corresponding
statement is usually false for the average values $\bar E_W$.

The {\em orbifold $E$-polynomial} can now be defined as
$$\Eorb (X/W) = \sum_{\{ w \} } \bar E_{C(w)}(X^w)(uv)^{F(w)}.$$
In the projective case, this reduces to 
$$\Eorb (X/W) = \sum_{p,q} h^{p,q}_{\oper{orb}}(X/W)(-u)^p(-v)^q.$$ 

The following theorem is proved by Batyrev-Dais \cite[6.14]{bd}.

\bs{Theorem}
\label{ind}
The orbifold $E$-polynomial depends only on the orbifold structure of
$X/W$.
\es

In fact, they show that it agrees with the more generally defined
notion of the ``stringy $E$-polynomial''.  Their proof as
stated assumes that $X$ is projective (or rather compact K\"ahler),
but it applies without change to the quasi-projective case.

\bit{Equality of the orbifold \boldmath $E$-polynomials}

We are ready to state and prove our main result, and then to show how
it implies the equalities predicted by mirror symmetry.

\bs{Theorem}
\label{mt}
Let $A$ be a connected abelian algebraic group, let $W$ be a finite
group acting orthogonally on a lattice $\La$ equipped with a positive
definite rational quadratic form, and let $B$ be a smooth variety
also acted on by $W$.  If $X
= (A \otimes \La) \times B$ and $\hat X = (A
\otimes \hat \La) \times B$ have the natural actions of
$W$, then $\Eorb(X/W) = \Eorb(\hat X/W)$.
\es

Even in the case $B=0$, the two quotients $X/W$ and $\hat X/W$ will
generally not be the same: see Proposition \re{exam}.  Moreover, the
theorem is certainly false if $\hat \La$ is replaced by some other
lattice isogenous to $\La$: see for example Theorem \re{sln}.

To prove the theorem, we will show that the sums defining the two
orbifold $E$-polynomials agree term by term: that is, for each $w \in
W$, the contributions of $X^w$ and $\hat X^w$ to the sums will be
identical.  To begin with, consider the fermionic shifts.

\bs{Lemma}
\label{r}
The fermionic shifts $F(w)$ and $\hat F(w)$ for the action of $w \in
W$ on $A \otimes \La$ and $A \otimes \hat \La$ are constant and equal.
\es

\pf.  Every fixed point $(A \otimes \La)^w$ is a subgroup.
Translation by any element of this subgroup commutes with the action
of $w$, so $F(w)$ is constant.  The Lie algebra of this subgroup is
$(\a \otimes \La)^w \subset (\a \otimes \La)$.  Since $\a$ is a
complex vector space, this equals $(\a \otimes_\C \h)^w$ where $\h =
\C \otimes \La$ is the Cartan subalgebra of $\g$.  The Killing form
induces an isomorphism $\h \cong \hat \h$; since $W$ preserves the
Killing form, it acts identically on these two spaces.  \fp

\bs{Lemma}
The identity components of $(A \otimes \La)^w$
and $(A \otimes \hat \La)^w$ have the same $E_{C(w)}$-polynomials.
\es

\pf.  Since the action of $w-1$ on $\La$ can be row-reduced over $\La
\otimes {\mathbb Q}$, a basis for $\ker(w~-~1) = \h^w \subset \h = \C
\otimes \La$ can be found in $\La$ itself.  Hence $\h^w = \C \otimes
\La^w$.  Consequently $(\a \otimes \La)^w = \a \otimes \La^w$.  The
identity component of $(A \otimes \La)^w$ is the connected subgroup of
$A \otimes \La$ with Lie algebra $(\a \otimes \La)^w$.  It therefore
is nothing but $A \otimes \La^w$.

Both $\La$ and $\hat \La$ can be regarded as subgroups of $\t = \R
\otimes \La$.  Since $\La^w$ and $\hat \La^w$ then have finite index
in the subgroup of $\t$ they jointly generate, both $A \otimes \La^w$
and $A \otimes \hat \La^w$ are quotients of the same abelian group by
a finite, $C(w)$-invariant subgroup.  This induces an isomorphism of
their compactly supported cohomology preserving the $C(w)$-action and
the mixed Hodge structure.  \fp

Denote the identity component of $(A \otimes \La)^w$ by $(A \otimes
\La)^w_0$.  The complex cohomology of $(A \otimes \La)^w$, as
a representation of $C(w)$ carrying a mixed Hodge structure, satisfies
$$H^*\lp(A \otimes \La)^w\rp \cong H^*\lp(A \otimes
\La)^w_0\rp \otimes \C[\pi_0((A \otimes \La)^w)].$$
Here $\pi_0$ denotes the group of components, and $\C[\phantom{x}]$
denotes a group algebra.
The equivariant $E$-polynomials therefore satisfy
$$E_{C(w)}\lp(A \otimes \La)^w\rp = E_{C(w)}\lp(A \otimes \La)^w_0\rp 
\,\, \chi\lp\C[\pi_0((A \otimes \La)^w)]\rp,$$
where $\chi$ denotes a character.

Theorem \re{mt} follows easily from this, together with the two lemmas
above, and the following.

\bs{Proposition}
\label{mp}
As representations of $C(w)$,
\end{em}
$$\C[\pi_0((A \otimes
\La)^w)] 
\cong
\C[\pi_0((A \otimes
\hat\La)^w)].$$

In fact we will see that the groups of components are canonically {\em
dual\/} to each other: that is, one is the group of characters
of the other, and the $C(w)$-action respects this duality.  

To determine the group of components $\pi_0((A \otimes \La)^w)$, it is
convenient to work over the real numbers.  Write $A$ as a product of
real lines and circles:
$A = \R^c \times \U{1}^d$, for some integers $c,d \geq 0$.  The
$W$-action on $A \otimes \La$ of course respects this decomposition, so
$$(A \otimes \La)^w 
= \lp(\R \otimes \La)^w\rp^c \times \lp(\U{1} \otimes \La)^w\rp^d.$$
The first factor is a linear subspace, so contributes nothing to the
component group.  As for the second factor, note that $\U{1} \otimes
\La$ is nothing but a maximal torus $T$ of the compact group $K$.  
So as sets with $C(w)$-action, \smallskip
\beq
\label{power}
\pi_0\lp(A\otimes \La)^w\rp = \pi_0\lp(T^w)^d\rp.
\eeq
We will not mention the action of $C(w)$ in what follows, but all of our
constructions are sufficiently natural to work
$C(w)$-equivariantly.

Let $\t = \R \otimes \La$.  Then there is an exact sequence of
$W$-modules 
$$0 \lrow \La \lrow \t \lrow T \lrow 0.$$
Moreover, the Killing form provides an inner product on $\t$.
An element $w \in W$ acts on $\t$ by a linear transformation $L$ which
is orthogonal with respect to the Killing form, so its action on the dual
$\hat \t$, via the adjoint inverse, is the same once we identify
$\t = \hat \t$.  On the other hand, there is also an exact sequence
$$0 \lrow \hat \La \lrow \hat \t \lrow \hat T \lrow 0.$$

Define $F: \t \to \t$ by $F(v) = L(v) - v$.  Then $T^w =
F^{-1}(\La)/\La$.  Similarly $\hat T^w = F^{-1}(\hat \La)/\hat\La$.  
Applying $F$ induces a natural homomorphism $\psi: T^w \to
\La/F(\La)$.

\bs{Lemma}
The map $\psi$ induces a natural isomorphism $\pi_0(T^w) =
\Tor(\La/F(\La))$.  
\es

\pf.  It suffices to show that $\ker \psi = T^w_0$ and $\oper{im} \psi =
\Tor(\La/F(\La))$.

Choose $v \in F^{-1}(\La)$ representing an element of $T^w =
F^{-1}(\La)/\La$.  If $\psi(v) = 0$, then $F(v) \in F(\La)$, so $v \in
\ker F + \La$.  Hence the element it represents in $F^{-1}(\La)/\La$
belongs to $\ker F / (\La \cap \ker F)$, which is precisely the
identity component of $F^{-1}(\La)/\La$.  Therefore $\ker \psi \subset
T^w_0$.  The other inclusion is easy.

Since $\psi$ factors through the finite group $T^w/\pi_0(T^w)$, its
image is certainly contained in the torsion part of $\La/F(\La)$.  On
the other hand, $\Tor(\La/F(\La)) = (\La \cap F(\t))/F(\La)$, since
$\Tor \La = 0$ and $F(\t)$ is the linear span of $F(\La)$.  And if $u
\in \La \cap F(\t)$, then certainly there exists $v \in F^{-1}(\La)$
such that $u = F(v)$.  Hence $\psi$ surjects onto
$\Tor(\La/F(\La))$. \fp 

\bs{Lemma}
There is a natural isomorphism \end{em}
$$\frac{\La}{F(\La)} 
= \Hom\left(\frac{F^{-1}_{\phantom{1}}(\hat
    \La)}{\hat \La},\U{1} \right).$$

\pf.  We show first that as subsets of $\t$, $F(\La) =
\Hom(F^{-1}(\hat \La), \Z)$.  For $x$ is in the right-hand side if and
only if $\langle x,y \rangle \in \Z$ whenever $y \in F^{-1}(\hat
\La)$, that is, whenever $\langle Fy,z \rangle \in \Z$ for all $z \in
\La$.  But
$$\langle Fy,z \rangle = \langle Ly,z \rangle - \langle y,z \rangle =
\langle y,L^{-1}z \rangle - \langle y,z \rangle = \langle y,F'z
\rangle$$
for $F' = L^{-1} - 1$.  So the right-hand side is double-dual to, and
hence equals, $F'(\La)$.  But $F' = -FL^{-1}$, and $L^{-1}$ preserves
$\La$, so $F'(\La) = F(\La)$. 

Now apply $\Hom(\phantom{x},\Z)$ to 
$$0 \lrow \hat \La \lrow F^{-1}(\hat \La) \lrow \hat T^w \lrow 0$$
to find
$$0 \lrow \Hom(\hat T^w,\Z) \lrow \Hom(F^{-1}(\hat \La),\Z) 
\lrow \Hom(\hat \La, \Z) \lrow \Ext^1(\hat T^w, \Z) \lrow 0.$$
Now $\Hom(\hat T^w,\Z) = 0$ since $\hat T^w$ is a compact abelian
group, and the next two terms have been identified as $F(\La)$ and
$\La$ respectively, so there is a natural isomorphism 
$$\Ext^1(\hat T^w, \Z) = \frac{\La}{F(\La)}.$$

On the other hand, apply $\Hom(\hat T^w, \phantom{x})$ to 
$$0 \lrow \Z \lrow \R \lrow \U{1} \lrow 0$$ to find
$$0 \lrow \Hom(\hat T^w, \Z) \lrow \Hom(\hat T^w, \R) \lrow \Hom(\hat
T^w, \U{1}) \lrow \Ext^1(\hat T^w, \Z) \lrow 0.$$
The first two terms are certainly 0 since $T^w$ is a compact abelian
group, so the last two terms are naturally isomorphic.  Applying
$\Hom(\phantom{x}, \U{1})$ to both sides yields the desired result.  
\fp

\bs{Lemma}
If $\Gamma$ is any finitely generated abelian group, such as
$\La/F(\La)$, then there is a natural isomorphism \end{em}
$$\pi_0(\Hom(\Gamma, \U{1})) = \Hom(\Tor \Gamma, \U{1}).$$

\pf.  The map is just restriction of homomorphisms to the torsion
part; that this is well-defined and induces the desired isomorphism
follows easily from the classification of finitely generated abelian
groups.  \fp

\noindent {\em Proof of Proposition \re{mp}}. 
Putting together the three preceding lemmas, we conclude that there is a
isomorphism 
$$\pi_0(\hat T^w) \cong \Hom(\pi_0(T^w),\U{1}),$$
which is natural, and, in particular, compatible with the $C(w)$-action.
Consequently, 
$$\C[\pi_0(\hat T^w)] \cong \C[\pi_0(T^w)]$$
as representations of $C(w)$: the map is given by the finite Fourier
transform $f \mapsto \hat f$, namely
$$\hat f(h) = \sum_{k \in \pi_0(T^w)} h(k) \, f(k).$$
Now take $d$th powers and apply \re{power}.  This completes the proof of
Proposition \re{mp}, and hence of Theorem \re{mt}. \fp

It is now easy to verify the predictions of mirror symmetry for the
normalized spaces of \S2.  Just take $A = \D$, $B=0$ in
Theorem \re{mt}.  This gives $\Eorb(\tilde M_\DR(G)) = \Eorb(\tilde
M_\DR(\hat G))$, as desired.  Similar identities hold for the
Dolbeault and Betti spaces.

If instead $A$ is taken to be an elliptic curve $D$, and $B$ to be $C
\otimes \La$ for another elliptic curve $C$, the theorem
shows that the spaces in the second row of the diagram from \S3 all
have the same $\Eorb$.  To draw the same conclusion
for the spaces in the first row, as predicted by mirror symmetry, we
need one last lemma, which shows that the hyperk\"ahler
rotation will not change $\Eorb$.

\bs{Lemma}
Suppose the Weyl group $W$ acts on a real torus $\mathcal B$ whose Lie
algebra $\mathfrak b$ is $W$-equivariantly isomorphic to $\Q \otimes
\La$, where $\Q$ are the quaternions.  If $\mathcal B$ is given the
hyperk\"ahler structure induced by $\Q$, then $\Eorb({\mathcal B}/W)$ is
the same for all complex structures in the hyperk\"ahler family.
\es

\pf.  All of the spaces involved in the computation of the
$E$-polynomials and fermionic shifts are tensor products of real
representations of $C(w)$ with the quaternions.  As complex
representations, they are therefore isomorphic in all complex
structures. \fp

\bit{Evaluation for quotients of \boldmath $\SL{n}$}

The orbifold $E$-polynomial $\Eorb((A \otimes \La)/W)$ described above
can in fact be evaluated more or less explicitly for the classical
groups.  For $\SL{n}$ it is computed by G\"ottsche \cite{g} and shown
to coincide with the Hodge polynomial of the crepant resolution,
namely the fiber of the natural sum morphism $\oper{Hilb}^n A \to A$.
For $\gp{Sp}{n}$ it is computed by Bryan-Donagi-Leung \cite{bdl} and
again shown to agree with the crepant resolution, the Hilbert scheme
of points on the K3 surface obtained by resolving $A/\pm 1$.  (These
papers assume $A$ is an abelian surface, but their computations work
for any 2-dimensional connected abelian algebraic group.)  One can
also obtain formulas for the Spin groups and for $\gp{SO}{2n}$, but
these are more cumbersome.

We present here a computation, based on that of G\"ottsche \cite{g}
and G\"ottsche-Soergel \cite{gs}, for $\Eorb((A \otimes\La)/W)$, where
$G$ is any quotient of $\SL{n}$, and $A$ is any connected abelian
algebraic group.  
  
So let $l$ and $m$ be positive integers with $lm=n$, let $G = \SL{n} /
\Z_m$, and let $\La$ be the coweight lattice.  The dual group is then
$\hat G = \SL{n} / \Z_l$, and the Weyl group is the symmetric
group $S_n$.  The conjugacy class of any $\si \in S_n$
is determined by the partition $\al$ of $n$ consisting of the lengths
of its cycles: say $n = \sum_j i_j$.  Let $\al_i$ be the number of
cycles of length $i$, so that $n = \sum_i i \al_i$ as well.  Let
$g=g(\al)$ be the greatest common divisor of those $i$ for which $\al_i
\neq 0$, and let $|\al| = \sum \al_i$ be the total number of cycles.
Recall also that $d$ is the number of $\U{1}$ factors
appearing in the decomposition $A = \R^c \times \U{1}^d$.

\bs{Theorem} 
\label{sln}
Let $A$, $\La$, and $W$ be as above.
Then as a polynomial in $u$ and $v$,
$$\Eorb\Left(\frac{A \otimes \La}{S_n} \Right) = 
\frac{1}{E(A)}\sum_{\al \in P(n)}
\tau_{l,m}^{g(\al),d}(uv)^{n-|\al|}  
\prod_i E(\oper{Sym}^{\al_i} A),$$
where $P(n)$ is the set of all partitions $n = \sum_i i \al_i$, and
$$\tau_{l,m}^{g,d} = \# \{ (r,s) \in \Z_g^d \times \hat \Z_g^d \st
(m,g)r=0=(l,g)s, \,\,
\langle r, s \rangle = 1 \}.$$
\end{em}

The $E$-polynomials $E(\oper{Sym}^a A)$ can easily be computed as 
$\bar E_{S_a}(A^a)$.  See G\"ottsche-Soergel \cite{gs} for a
convenient formula when $A$ is an abelian surface. \medskip

\pf. Identify the coweight lattice of $\GL{n}$, as an $S_n$-module, with
$\Z^n$.  Then the coweight lattice of $\SL{n}$ is the kernel of the
sum map $\Z^n \to \Z$; the coweight lattice of $\GL{n}/\Z_m$ is the
lattice generated by $\Z^n$ and the point whose coordinates are all
$1/m$; and the coweight lattice of $G$, namely $\La$, is the kernel of
the latter lattice under the sum map.

Therefore, if $K$ is the kernel of the sum map $A^n \to A$, then as an
$S_n$-module, $A \otimes \La = K/A[m]$, where $A[m] \cong \Z_m^d$
denotes the division points of order $m$ in $A$, acting diagonally on
$A^n$ and hence on $K$.  The quotient $(A \otimes \La)/S_n$ can then
be identified with $K/(A[m] \times S_n)$.  By Theorem \re{ind}, the
orbifold $E$-polynomial can equally well be computed from this
quotient.

So fix $\si \in S_n$ and $a \in A[m]$.  
Also, let $k$ be the order of $a$ in $A[m]$, so that $k|m$.  To be
fixed by the action of $\si \times a$, an $n$-tuple in $A^n$ must
consist of cycles of the form 
\beq \label{cycle}
(x+a,x+2a,x+3a,\dots,x+ia) 
\eeq
with $x+ia=x$, whose number and size are given by $\al$.  Hence,
for $(A^n)^{\si \times a}$ to be nonempty, it is necessary and
sufficient that $k|i$ for every $i$ such that $\al_i \neq 0$, and
hence that $k|g$.  Furthermore, if this is satisfied, then $(A^n)^{\si
\times a}$ can be identified with $A^{|\al|}$ by choosing the first
element of each cycle.

However, the intersection of $(A^n)^{\si \times a}$ with $K$ may not be
connected.  Indeed, the sum map $A^n \to A$ restricts to
$(A^n)^{\si \times a} = A^{|\al|}$ as $(x_j) \mapsto qa + gf_\al(x_j)$,
where $f_\al: A^{|\al|} \to A$ is
$$f_\al(x_j) = \sum_j \frac{i_j}{g}\, x_j,$$
and $q$ is a constant,
either $0$ or $k/2$, depending on a straightforward parity condition.
Consequently, $(x_j) \in (A^n)^{\si \times a}$ belongs to $K$ if and
only if $f_\al(x_j)$ is a $g$th root of $-qa$.  Such $g$th roots, of
course, form a coset of $A[g]$, which partitions $K^{\si \times a}$
into $g^d$ disjoint subvarieties, each isomorphic to $K_\al = \ker
f_\al$.

These subvarieties are actually connected, and so constitute the
components of $K^{\si \times a}$.  Indeed, 
since the $i_j/g$ are coprime, there exist $c_j \in \Z$ such that
$\sum_j c_j i_j/g = 1$.  The map $A \to A^{|\al|}$ given by $x \mapsto
(c_j x)$ is then a right inverse for $f_\al$ (cf.\ \cite{gs}), so that
$K_\al \times A \cong A^{|\al|}$.

In fact, more is true.  Let $\Ga$ be the subgroup $\prod_i S_{\al_i}
\subset C(\si)$ interchanging the cycles of the same length.  Then
$\Ga$ acts on $A^{|\al|}$ preserving $K_\al$, and the action commutes
with the isomorphism above, if the second factor of $K_\al \times A$
is given the trivial $\Ga$-action.  Hence 
\beq \label{prod}
E_\Ga(K_\al) = \frac{E_\Ga(A^{|\al|})}{E(A)}. \eeq

The action of the full centralizer $C(\si \times a) = C(\si) \times
A[m]$ does not preserve the components $K_\al$ of $K^{\si \times a}$,
however.  First consider the action of $A[m]$.  The addition of $b \in
A[m]$ to each $x_j$ changes $f_\al(x_j)$ to $f_\al(x_j+b)$ =
$f_\al(x_j) + (n/g)b$.  So the components of $K^{\si \times a}/A[m]$
are indexed by $A[g]/(n/g)A[m] \cong A[(l,g)]$, where $l=n/m$.

As for the action of $C(\si)$, it is convenient to break it into two
pieces.  Let $\Si \subset C(\si)$ be the normal subgroup $\prod_j
\Z_{i_j}$ of $C(\si)$, acting cyclically on the cycles of $\si$.  Note
that $C(\si)/\Si$ is the group $\Ga$ mentioned before.  Since every
element of $K^{\si \times a}$ has cycles of the form shown in
\re{cycle}, the generator of $\Z_{i_j} \subset \Si$ acts by adding $a$
to $x_j$.  Since the $i_j/g$ are coprime, the action of $\Si$ can thus
add any multiple of $a$ to $f_\al(x_j)$.  Hence the components of
$K^{\si \times a}/(A[m] \times \Si)$ are indexed by the quotient
$A[g]/((n/g)A[m] + \langle a \rangle)$.

Counting the number of elements in this finite abelian group is
facilitated by passing to the dual: $A[g] \cong \Z_g^d$ is dual to
$\hat A[g] = \Hom(A[g], \U{1})$, and so this quotient is dual to the
subgroup of $\hat A[(l,g)] \subset \hat A[g]$ evaluating to $1$ on
$a$.  Summing over $a \in A[(m,g)]$, we find that the total number of
components of
$$\bigcup_{a \in A[(m,g)]}\frac{K^{\si \times a}}{A[m] \times \Si}$$
is given by the number $\tau_{l,m}^{g,d}$ defined in the statement.
Clearly $\tau_{l,m}^{g,d}$ is symmetric in $l$ and $m$, in agreement
with Theorem \re{mt}.  When $(l,g) = 1$ it equals $(m,g)^d$, which is
then just $g^d$ since $g | lm$.

The part of $\Si \times A[m]$ which preserves $K_\al$ acts merely by
translations, so dividing by it has no effect on the $E$-polynomial.
Hence
$$\bar E_{C(\si) \times A[m]}
\Left( \bigcup_a K^{\si \times a} \Right) = 
\tau_{l,m}^{g(\al),d} \bar
E_\Ga(K_\al).$$
The action of $\Ga$, on the other hand, permutes the cycles having the
same length.  The quotient $A^{|\al|} / \Ga$ is therefore $\prod_i
\oper{Sym}^{\al_i} A$, so according to \re{prod},
$$\bar E_\Ga(K_\al) = \frac{\prod_i E(\oper{Sym}^{\al_i} A)}{E(A)}.$$

It remains only to compute the fermionic shifts.  But these are easily
seen to be independent of $a$, and in the case $a=0$ are computed by
G\"ottsche \cite{g} and G\"ottsche-Soergel \cite{gs} to be just
$n-|\al|$. \fp

\bit{Some final remarks} 

\vspace{-12pt}
\br{Analogy with the Greene-Plesser construction} 
In the case where $G$ is semisimple, simply connected and simply
laced, that is, of type $A$, $D$, or $E$, then the Langlands dual
$\hat G$ is simply $G / Z(G)$, and $\hat \H = \H/Z(G)$.  Consequently,
the spaces appearing in the second columns of the diagrams of \S\S2 and
3 --- that is, the mirrors of the first columns --- are quotients of
those in the first columns by the finite abelian group $Z(G)^d$.  This
is eerily reminiscent of the Greene-Plesser construction \cite{gp} of
the mirror of the Fermat quintic $Q$ as a quotient of $Q$ by the
finite abelian group $\Z_5^4$.  Perhaps this analogy could be expanded
to include other groups.  Conceivably these spaces could even be realized
in a suitable way as subvarieties of dual toric varieties.  \er

\br{Topologically non-trivial bundles}
This paper deals with only one component of the moduli spaces it
employs.  An obvious question is how much carries over to other
components.  Further, one could ask the same question about
spaces of bundles which are topologically non-trivial.  This would
undoubtedly be more difficult.  In particular, the work of
Borel-Friedman-Morgan \cite{bfm} would come into play. \er

\br{Hilbert schemes of points}  When $G=\SL{n}$, the quotient $(A
\otimes \La)/W$ has a crepant resolution, namely the identity fiber
$K_{n-1}A$ of
the natural map $\oper{Hilb}^n A \to A$.  
It does not follow from the results of this paper, but it seems
extremely plausible, that for $m|n$ 
$$\Eorb(K_{n-1}A/\Z_m^d) 
= \Eorb((A \otimes \La)/(W \times \Z_m^d)).$$
After all, the left-hand side is a partial crepant resolution of the
right-hand side.  This should be straightforward to check, following
G\"ottsche \cite{g}. 

\br{Equivalence of derived categories} 
In 1994, Kontsevich proposed an extension of mirror symmetry, the
so-called ``homological mirror symmetry'' \cite{k}.  For a mirror
pair $M, \hat M$ of Calabi-Yau 3-folds, it is supposed to identify two
derived categories: that of coherent sheaves on $M$, and that of the
Fukaya category on $\hat M$.  In even dimensions, we expect mirror
symmetry to identify the A-model of $M$ with the A-model, not the
B-model, of $\hat M$.  Hence for our examples, we would expect the
derived categories of coherent sheaves to be isomorphic for $M(G)$ and
$M(\hat G)$, and likewise for the Fukaya categories.

Even the definition of the Fukaya category remains somewhat
mysterious, so there is no hope of proving this equivalence.  
One might wonder whether hyperk\"ahler rotation transforms the Fukaya
category into the category of coherent sheaves, but this is pure
speculation. 

However, one equivalence, that of derived categories of coherent
orbifold sheaves on $\tilde M_\Dol(G)$ and $\tilde M_\Dol(\hat G)$, is
easy.  Or rather, it follows immediately from the work of Mukai
\cite{m1}.  He showed that the derived categories of coherent sheaves
on an abelian variety and its dual are equivalent.  The equivalence is
defined by pulling back to the product, tensoring by the Poincar\'e
line bundle and pushing forward to the other factor.  This can be done
for $C \otimes \La$ over the base $\h$, and it can be done
$W$-equivariantly, inducing an equivalence of derived categories for
sheaves with $W$-action.  The same works for the spaces of \S3.  Note,
however, that even for $G=\SL{2}$ the Poincar\'e line bundle does not
descend to a bona fide sheaf on the quotient by $W$.  Using orbifold
sheaves therefore seems necessary.  \er

\end{document}